\documentclass[11pt]{article}

\usepackage{amsmath}
\usepackage{graphics}
\usepackage{epsfig}
\usepackage{amssymb}
\usepackage{amsfonts}
\usepackage{color}

\newtheorem{lem}{Lemma}
\newtheorem{teo}{Theorem} 

\newcommand{\PT}{\mathcal{PT}}

\begin{document}
\title{Supercritical blowup in coupled  parity-time-symmetric nonlinear Schr\"{o}dinger equations}

\author{Jo\~ao-Paulo Dias$^{1}$, M\'ario Figueira$^{1}$,  Vladimir V. Konotop$^{2,3}$, and Dmitry A. Zezyulin$^{2}$
\\
$^{1}$ CMAF/UL, Faculdade de Ci\^encias,
Universidade de Lisboa,
\\
Avenida Professor
Gama Pinto 2, Lisboa 1649-003, Portugal
\\
$^{2}$
 Centro de F\'isica Te\'{o}rica e Computacional,
 Faculdade de Ci\^encias, Universidade de Lisboa,
\\
Avenida Professor
Gama Pinto 2, Lisboa 1649-003, Portugal
\\
$^{3}$
 Departamento de F\'isica,  Faculdade de Ci\^encias, Universidade de Lisboa,
\\
Campo Grande, Ed. C8,
Lisboa 1749-016,  Portugal
}

\maketitle

\begin{abstract}
We prove finite time supercritical blowup in a parity-time-symmetric system of the  two coupled nonlinear Schr\"odinger  (NLS) equations.  One of the equations contains gain and the other one contains  dissipation such that strengths of the gain and dissipation are  equal.  We address two  cases:  in the first model all nonlinear coefficients  (i.e. the ones describing self-action and non-linear coupling) correspond to attractive (focusing) nonlinearities, and  in the second case the NLS equation with gain has attractive nonlinearity  while the NLS equation with dissipation has repulsive (defocusing) nonlinearity and the nonlinear coupling is repulsive, as well. The proofs are based on the virial technique arguments. Several   particular cases are also illustrated numerically.
\end{abstract}

\vspace{0.5cm}

\section{Introduction and preliminary observations}

In this paper we consider the Cauchy problem for two coupled nonlinear Schr\"odinger (NLS) equations
\begin{subequations}
\label{Cauchy}
\begin{eqnarray}
\label{NLS1}
&i u_t =-\Delta u +i\gamma u+\kappa v-(g_1|u|^2+g|v|^2)u,
&
\\[1mm]
\label{NLS2}
&i v_t=-\Delta v -i\gamma v+\kappa u-(g|u|^2+g_2|v|^2)v, &
\end{eqnarray}
where   $x\in\mathbb{R}^N$, $t\geq 0$, $\Delta$ is the $N$-dimensional Laplacian, and
\begin{eqnarray}
\label{IC1}
\begin{array}{c}
u(x,0)=u_0(x), \quad v(x,0)=v_0(x), \quad
u_0,\, v_0\in H^{1}(\mathbb{R}^N),
\\[2mm] |x|\,u_0(x),\,\, |x|\,v_0(x)\in L^2(\mathbb{R}^N).
\end{array}
\end{eqnarray}
\end{subequations}

We  consider
 $g_1,g_2,g\in\mathbb{R}$,
and assume that $\gamma > 0$ and $\kappa>0$ [the cases $\gamma<0$ and (or) $\kappa<0$ can  be considered analogously after simple redefinition of the functions $u(x,t)$ and $v(x,t)$].
%
%
Our main interest is in establishing  sufficient conditions for the initial conditions $u_0(x)$ and $v_0(x)$ to blow up in  finite time.

The system (\ref{NLS1})--(\ref{NLS2}) is  referred to as parity-time ($\PT$) symmetric because of its linear counterpart ($g=g_1=g_2=0$), which under certain conditions allows for stable propagating of the linear waves, which  justifies the particular relevance of the model. The concept of $\PT$ symmetry can be formulated also in the nonlinear case referring to the formal property as follows: if functions $u(x,t)$ and $v(x,t)$ solve the  pair of equations in   (\ref{NLS1})--(\ref{NLS2}), then the new functions $u_{\PT}(x,t):=\bar{v}(x,-t)$ and $v_{\PT}(x,t):=\bar{u}(x,-t)$ also solve the same equations, provided that $g_1=g_2$ (hereafter an overbar stands for the complex conjugation).

The concept of  $\PT$ symmetry was originally employed in quantum mechanics to construct non-Hermitian potentials
with pure real spectra~\cite{Bender,PT}.
Soon after that, the ideas of  $\PT$ symmetry were transferred to the   optics~\cite{Muga}. Using the  mathematical analogy between the Schr\"{o}dinger equations in quantum mechanics and paraxial approximation in optics, it was shown that a medium with mutually balanced gain and losses  allows for   stable propagation of linear waves. Further, it was suggested in~\cite{OL} that a simple implementation of $\PT$ symmetry in optical systems can be achieved using two waveguides. This geometry was proven to be very suitable   for experimental observation of the effects related to the $\PT$ symmetry~\cite{exp0,exp1}. Optical applications which naturally admit inclusion of the nonlinearity have stimulated further studies of the nonlinear $\PT$-symmetric system~\cite{Musslimani}. In  the two-waveguide geometry described by the coupled NLS equations with gain and dissipation,  there have been considered bright~\cite{bright1,bright2} and dark~\cite{dark} solitons, breathers~\cite{breather}, and rogue waves~\cite{rogue}.

Interplay between the $\PT$ symmetry and nonlinearity  was also considered in the context of Bose--Einstein condensates (BECs)~\cite{Cartarius,Cartarius2}, where the coupled NLS (alias Gross--Pitaevskii) equations appear as a natural model for two BECs split in an external double-well trap. Alternatively, the $\PT$ symmetry in BECs can be implemented using spinor components of a single atomic specie in two different  ground states, one of which is pumped by an external laser and another one has (either natural or induced) loss of atoms~\cite{SO-BEC}. Since a three-dimensional realization of a BEC appears in the most of experimental settings, it is of interest to consider dynamics of two coupled three-dimensional, and more generally $N$-dimensional, NLS equations (\ref{Cauchy}) one of which accounts for  gain [equation~(\ref{NLS1})] and the other one has dissipation [equation~(\ref{NLS2})].

Turning now to the blowup phenomenon, in coupled conservative NLS equations it was previously discussed in Refs.~\cite{coupled1,coupled2} without linear coupling  and in Ref.~\cite{coupled3} with linear coupling (see also~\cite{Berge} and references there in for the critical collapse). On the other hand, blowup in a single NLS equation with homogeneous and inhomogeneous linear damping was considered in \cite{OT,Tsutsumi} and  \cite{DF}, respectively.

As mentioned above, particular interest in analyzing the blowup in the $\PT$-symmetric  coupled NLS equations (\ref{NLS1})--(\ref{NLS2}) stems from the fact that depending on the relation among parameters the system can display features either of a conservative medium, where linear waves [i.e. solutions of the Eqs.~(\ref{Cauchy}) with $g_1=g_2=g=0$] can propagate stably or of a dissipative one where a linear mode is unstable (i.e. it either decays or infinitely grows with time). The former case is referred to as the {\em unbroken $\PT$-symmetric phase} \cite{Bender}, and corresponds to
\begin{eqnarray}
\label{unbroken}
\quad \kappa> \gamma >0.
\end{eqnarray}
If $\gamma >\kappa>0$, then the  $\PT$ symmetry   is said to be broken. The particular case $\kappa=\gamma$ corresponds to the exceptional point~\cite{Kato} of the underlying linear operator.

Thanks to the possibility of stable propagation of linear modes, the case of unbroken phase is of special physical interest. It also allows one to make some general preliminary conclusions about blowing up and global existence of solutions. Indeed, if (\ref{unbroken}) is satisfied one can perform the transformation
\begin{eqnarray}
\label{rotation}
\left( \begin{array}{c}
u\\v
\end{array}\right)
=\left( \begin{array}{cc}
e^{i\alpha/2} & -e^{-i\alpha/2}\\ e^{-i\alpha/2} & e^{i\alpha/2}
\end{array}\right)\left( \begin{array}{c}
U\\V
\end{array}\right),
\end{eqnarray}
where the constant phase $\alpha$ is determined by the relation
\begin{eqnarray}
e^{i\alpha}=-\frac{\kappa}{\sqrt{\kappa^2-\gamma^2}-i\gamma}.
\end{eqnarray}
We notice that (\ref{unbroken}) implies that $\alpha\in [0,\pi/2)$.
The newly introduced functions  $U(t)$ and $V(t)$ solve the system
 \begin{eqnarray}
\label{NLS1_rot}
i U_t=-\Delta U +\omega U-F_1,
\quad
i V_t=-\Delta V -\omega V-F_2
\end{eqnarray}
with $\omega=\sqrt{\kappa^2-\gamma^2}$,
\begin{eqnarray*}
F_1=G|U|^2U+G_+|V|^2U+2Q|U|^2V+M|V|^2V+PU^2\bar{V}+(G-2g)V^2\bar{U},
\\
F_2=\bar{G}|V|^2V+G_+|U|^2V+2\bar{Q}|V|^2U+M|U|^2U+\bar{P}V^2\bar{U}+(\bar{G}-2g)U^2\bar{V},
\end{eqnarray*}
and the other parameters being defined by
\begin{eqnarray*}
G_\pm=\frac 12 (g_1\pm g_2),
\quad
G=g+G_+ -iG_-\tan \alpha,
\\
M=-\frac{G_-}{\cos\alpha},
\quad
Q=-M+ig\sin \alpha,
\quad
P=M\cos(2\alpha)-2iG_+\sin\alpha
\end{eqnarray*}
and the initial conditions
\begin{eqnarray}
U(x,0)=U_0(x), \quad V(x,0)=V_0(x),\quad U_0,\, V_0\in H^{1}(\mathbb{R}^N).
\end{eqnarray}
Notice that the performed transformation (\ref{rotation}) does not depend on the dimensionality of the space $N$.

If $g_1=g_2$ we have  $G_-=M=0$ and the corresponding system (\ref{NLS1_rot})
has  particular solutions which can be found in the form $V\equiv 0$  with $U$ solving the equation
\begin{eqnarray}
\label{NLS_cut}
i U_t=-\Delta U +\omega U- (g+g_1)|U|^2U.
\end{eqnarray}
(analogously one can consider a particular solution with    $U\equiv 0$). This leads to a number of conclusions (which follow from the well known results on the NLS equation, see. e.g.~\cite{Sulem}) as follows.

First, if $g+g_1>0$ and $N\geq 2$, then there exist solutions blowing up in finite time. This collapse is characterized by the conserved squared norms  $\|u\|_2^2=\|v\|_2^2=\|U\|_2^2$  and simultaneous blowup of  both the fields $u$ and $v$.  Notice that hereafter we use
the abbreviated notation for the standard $L^n$-norm: $\|\cdot\|_n= \|\cdot\|_{L^n}$.

Second, if $g+g_1<0$ then there exist  global (dispersive) solutions.

Third, subject to sufficiently small initial conditions and $N=2$ (the critical case)  one can find solutions existing globally.  Moreover, in the case $g_1=g_2=-g$ a solution with either $V\equiv 0$ or $U\equiv 0$ with smooth localized initial data exists globally, since it is described by the linear Schr\"odinger equation.

Below in this paper we concentrate on $N\geq 3$ corresponding to the super-critical case.

\section{Main results}
\label{sec:main}
In order to formulate our main results it is convenient to define the following integral Stokes components for the NLS equation:
\begin{eqnarray*}
S_0=\|u\|_2^2+\|v\|_2^2,\quad S_1=2\mbox{Re}\!\int \! u\bar{v}\,dx, \\
S_2=2\mbox{Im}\!\int\! u\bar{v}\,dx, \quad S_3=\|u\|_2^2-\|v\|_2^2.
\end{eqnarray*}
Hereafter for the sake of brevity  we use the notation $\int \ldots dx =\int_{\mathbb{R}^N} \ldots  dx$.
Obviously, $S_0$   is a conserved quantity for $\gamma=0$, but becomes  time-dependent for $\gamma>0$. Indeed, it is   straightforward to compute
\begin{eqnarray}
\label{difS0}
\frac 12 \frac{d\,S_0(t)}{dt}
=\gamma
S_3(t),
\end{eqnarray}
which readily gives the  estimate
\begin{eqnarray}
\label{S_estim}
S_0(t)\leq S_0(0) e^{2\gamma t}.
\end{eqnarray}

We also define the energy
\begin{eqnarray}
\label{energy}
E(t)=\int\left(|\nabla u|^2+|\nabla v|^2 +\kappa(u\bar{v}+\bar{u}v) -\frac{g_1}{2}|u|^4-\frac{g_2}{2}|v|^4-
-g|u|^2|v|^2\right)dx
\end{eqnarray}
which is another conserved quantity in the conservative case ($\gamma=0$). For arbitrary  $\gamma$   the energy changes with time   according to
\begin{eqnarray}
\label{DE}
\frac{dE}{dt}=2\gamma \int \left(|\nabla u|^2-|\nabla v|^2 - g_1|u|^4+g_2|v|^4\right)dx.
\end{eqnarray}
Next, we define
 the mean squared width of the solution $X$
\begin{eqnarray}
X(t)=\int |x|^2\left(|u|^2+|v|^2\right)dx,
\end{eqnarray}
and its derivative $Y$:
\begin{eqnarray*}
\label{Y}
Y(t)\equiv \frac{dX(t)}{dt}=4\mbox{Im}\int(ux\cdot\nabla\bar{u}+vx\cdot\nabla\bar{v})dx+2\gamma\int |x|^2(|u|^2-|v|^2)dx.
\end{eqnarray*}
For the second derivative of $X(t)$ one can compute
\begin{eqnarray}
\label{X} \frac{d^2 X(t)}{dt^2}=4NE(t)+4\int(2-N)\left(|\nabla
u|^2+|\nabla v|^2\right)dx
+4\gamma^2X+
\nonumber\\
+16\gamma\mbox{Im}\int\left(\bar{u}x\cdot\nabla
u-\bar{v}x\cdot\nabla
v\right)dx+8\gamma\kappa\mbox{Im}\int|x|^2v\bar{u}dx-
\nonumber\\
-8\kappa N\mbox{Re}\,\int u\bar{v}dx.
\end{eqnarray}

Now  we introduce   functions
\begin{eqnarray}
\label{F}
F(t) &=& X(0)+Y(0)t+\frac{8N}{N+2}E(0)t^2
+\frac{4\kappa}{\gamma^2}
S_0(0)
\left(e^{2\gamma t}- 2\gamma t - 1\right),\\[2mm]
\label{M}
M(t)& =& \sup_{\tau\in[0,t]}F(\tau)+1>0,
\\[2mm]
\label{G}
G(t)  &=&  M(t)\left(c_1 \frac{t^2}{2}+
\exp\left(\frac{c_3\gamma t}{c_2}\right) - 1\right),
\end{eqnarray}
where $S_0(0)=\|u_0\|_2^2+\|v_0\|_2^2 $ is defined by  $L^2(\mathbb{R}^N)$-norms of the initial conditions $u_0$ and $v_0$,
as well as   constants
\begin{eqnarray*}
c_1&=&4\gamma\kappa+4\gamma^2\frac{ {5N+6}}{N-2},\\[3mm]
c_2&=&\left\{
\begin{array}{cl}
\frac{4}{5}\min\{1, g_1, g_2\}, &\quad \mbox{if} \quad g\geq 0,\\[2mm]
\frac{4}{5}\min\left\{1, g_1+\frac{g\sqrt{g_1}}{\sqrt{g_2}},
g_2+\frac{g\sqrt{g_2}}{\sqrt{g_1}}\right\}, &\quad \mbox{if} \quad
-\sqrt{g_1g_2} < g < 0,
\end{array}
\right.
\\[3mm]
c_3&=&\frac{32 N }{N+2}\max\{1,g_1,g_2\}
\end{eqnarray*}

Our main result consists in the following  Theorem which provides sufficient conditions for the  finite-time blowup of solution for the problem (\ref{Cauchy}) in the supercritical case.
\begin{teo}
\label{theor1} Let $N\geq 3$ and
\begin{eqnarray}
\label{g_positive}
g_{1,2}>0, \qquad g>-\sqrt{g_1g_2}.
\end{eqnarray}
 Assume
that the initial conditions $u_0(x)$ and $v_0(x)$ of the Cauchy problem  (\ref{Cauchy}) are chosen such that there exists
$T_0>0$ for which the following two conditions hold:
\begin{eqnarray}
\label{cond1}
 F(T_0)+1<0,\\
 \label{cond2}
G(T_0) < 1.
\end{eqnarray}
Then the solution of the problem (\ref{Cauchy}) does not exist in the
interval $t\in [0,T_{0}]$.
\end{teo}

Before presenting the proof of the Theorem~\ref{theor1}, it is useful to ensure that the conditions (\ref{cond1}) and (\ref{cond2}) are consistent and that the   initial conditions $u_0$ and $v_0$ satisfying  (\ref{cond1})--(\ref{cond2}) do exist. We illustrate the consistency of the conditions (\ref{cond1}) and (\ref{cond2}) formulating  the two lemmas as follows. Lemma~1  establishes that the blow up occurs if the initial energy is negative (with large enough absolute value), while Lemma~2 shows that the blow up conditions can be satisfied if $Y(0)$ is   negative (with large enough absolute value).
\begin{lem}
\label{lem1} Let $\beta =\frac{c_3\gamma}{c_2}$,
\begin{equation*}
{C_0} = \frac{|Y(0)|}{2\gamma} + \frac{4\kappa  }{\gamma^2}
S_0(0), \quad \tilde{M}(t) = 1 + X(0)+{C_0}(e^{2\gamma t} - 1),
\end{equation*}
and
\begin{eqnarray}
{T}_{0,max} = \frac{1}{\beta} \ln\left (1 +  \frac{\beta^2}{(1+X(0))\left(\beta^2+ c_1 \right)}\right),\\
{T}_{0,min} = \frac{1}{\beta} \ln\left (1+  \frac{\beta^2}{\tilde{M}(T_{0,max})\left(\beta^2+c_1\right)}\right).
\end{eqnarray}
Then there exists  $T_0\in [ {T}_{0,min},{T}_{0,max}]$ such that the blow-up conditions (\ref{cond1})--(\ref{cond2}) are satisfied at $t=T_0$ provided that
\begin{equation}
\label{E0suff} E(0) < -
\frac{(N+2)\tilde{M}({{T}_{0,max}})}{8N{T}_{0,min}^2}.
\end{equation}
\end{lem}

\textbf{Proof.}
Let us introduce
\begin{eqnarray}
\tilde{G}(t) :=  {M}(t)\left(1 + \frac{c_1}{\beta^2}\right)\left(e^{\beta t} - 1\right) >  G(t) \mbox{\quad for } t>0,
\end{eqnarray}
and define
 ${T}_0$ as  the smallest solution of
  the equation  $\tilde{G}(T_0)=1$ (which implies that $G(T_0) < 1)$.
 Using that $M(t) \geq 1+X(0)$ we find that $T_0 \leq T_{0,max}$. Introduce $\tilde{F}(t) = F(t) - \frac{8N}{N+2}E(0) t^2$. For $E(0) \leq 0$ we have   $\tilde{F} \geq F(t)$.  With simple transformations we obtain that $\tilde{F}(t) + 1 \leq  \tilde{M}(t)$ and therefore for all $t\in [0, T_{0,max}]$ one has
 \begin{equation}
 \label{boundM}
 M(t) \leq \sup_{\tau \in [0, t]} \tilde{F}(\tau) + 1 \leq \sup_{\tau \in [0, t]} \tilde{M}(\tau)  = \tilde{M}(T_{0,max}),
 \end{equation}
 which implies that  $T_0 \geq T_{0,min}$.

In order to satisfy the condition (\ref{cond1}) we require
$0>F({T}_0)+1 = \tilde{F}({T}_0)+\frac{8N}{N+2}E(0) {T}_0^2 + 1$
which is satisfied automatically if
\begin{equation}
\label{E0suff1} E(0) < -
\frac{(N+2)\tilde{M}({{T}_{0}})}{8N{T}_{0}^2}.
\end{equation}
If the condition (\ref{E0suff}) holds then the   requirement (\ref{E0suff1}) also holds. Therefore, both conditions of the blow-up are satisfied at $t=T_0$.
\hfill $\blacksquare$

\begin{lem}
\label{lem2}
Redefine the constant $C_0$ as follows:
\begin{eqnarray}
{C_0} = \frac{4N|E(0)|}{(N+2)\gamma^2} + \frac{4\kappa}{\gamma^2}
S_0(0),
\end{eqnarray}
and keep others definitions the same as in Lemma~\ref{lem1}. The conditions (\ref{cond1})--(\ref{cond2})  are satisfied  if
\begin{equation}
Y(0) < \frac{8\kappa  }{\gamma}S_0(0)
- \frac{\tilde{M}({{T}_{0,max}})}{{T}_{0,min}}.
\end{equation}
\end{lem}
\textbf{Proof.} The proof is almost identical to that for Lemma~\ref{lem1} except for  definition of  the function $\tilde{F}(t)$. Now it is defined as $\tilde{F}(t)= F(t) - Y(0)t + 8\kappa   S_0(0) t
/ \gamma$. Then  $\tilde{F}(t) \geq F(t)$ for any $Y(0) \leq
\frac{8\kappa   }{\gamma}S_0(0)$. \hfill $\blacksquare$

\section{Proof of Theorem~\ref{theor1}}

We start with the following estimate which is obtained from (\ref{DE}) and from the definition of the constant $c_3$:
\begin{eqnarray}
\label{E_estim}
E(t)\leq 2\gamma
 \int_{0}^{t}\left(\|\nabla u\|_2^2+\|\nabla
v\|_2^2+g_1\|u\|_4^4+g_2\|v\|_4^4\right)d\tau+E(0)
\nonumber \\
\leq \gamma c_3
\frac{N+2}{16N}\int_{0}^{t}\left(\|\nabla u\|_2^2+\|\nabla
v\|_2^2+\|u\|_4^4+\|v\|_4^4\right)d\tau+E(0).
\end{eqnarray}
Rearranging the terms two in the r.h.s. of (\ref{X}) as follows
\begin{eqnarray}
\label{rearrange}
4NE(t)+4\int(2-N)\left(|\nabla u|^2+|\nabla v|^2\right)dx=
\frac{1}{N+2}\Biggl\{
{16N}E(t)-
\nonumber\\
{8(N-2)}\int\left(|\nabla u|^2+|\nabla v|^2\right)dx-
 2N{(N-2)}\left[g_1\int|u|^4dx + g_2\int |v|^4dx \right] -
 \nonumber\\
4Ng{(N-2)}\int|u|^2|v|^2dx +
 8N{(N-2)}\kappa\mbox{Re}\int u\bar{v}dx \Biggr\},
\end{eqnarray}
and taking into account that
\begin{eqnarray}
\label{auxil_1}
16\gamma\mbox{Im}\int\left(\bar{u}x\cdot\nabla u-\bar{v}x\cdot\nabla v\right)dx\leq
\frac{4(N-2)}{N+2}\left(\|\nabla u\|_2^2+\|\nabla v\|_2^2\right) +
\frac{16\gamma^2(N+2)}{N-2}X,
\end{eqnarray}
and
\begin{eqnarray}
8\gamma\kappa \mbox{Im}\int|x|^2v\bar{u}\,dx\leq 4\gamma\kappa X,
\end{eqnarray}
from (\ref{X})  we obtain the following estimate:
\begin{eqnarray}
\label{intermediate}
\frac{d^2 X}{dt^2} + \frac{2(N-2)}{N+2}\left(2\|\nabla
u\|_2^2+2\|\nabla v\|_2^2 + Ng_1\|u\|_4^4+Ng_2\|v\|_4^4+ 2Ng \int
|u|^2|v|^2dx \right) \leq
\nonumber\\
\leq c_1X+\frac{16N}{N+2}E(t)+{16 \kappa }S_0(t),
\end{eqnarray}
which implies
\begin{equation}
\frac{d^2 X}{dt^2} +c_2\left(\|\nabla u\|_2^2+\|\nabla v\|_2^2+\|u\|_4^4+\|v\|_4^4\right)\leq
 c_1 X+\frac{16N}{N+2}E(t)+ 16 \kappa  S_0(t).
\end{equation}
Then using (\ref{S_estim}) and (\ref{E_estim})  we arrive at
\begin{eqnarray}
\label{X2} \frac{d^2 X}{dt^2} +c_2\left(\|\nabla u\|_2^2+\|\nabla
v\|_2^2+\|u\|_4^4+\|v\|_4^4\right)\leq  c_1X+\frac{16N}{N+2}E(0)+
\nonumber \\
 + 16\kappa
S_0(0)
e^{2\gamma t} +  c_3 \gamma \int_{0}^{t}\left(\|\nabla
u\|_2^2+\|\nabla v\|_2^2+\|u\|_4^4+\|v\|_4^4\right)d\tau
\end{eqnarray}
Next, introducing
\begin{eqnarray}
\rho(t)=\int_{0}^{t}\int_{0}^{ \sigma}
\left(\|\nabla u\|_2^2+\|\nabla v\|_2^2+\|u\|_4^4+\|v\|_4^4\right)d\tau d\sigma,
\end{eqnarray}
and using function $F(t)$  defined in  (\ref{F}), 
we rewrite (\ref{X2}) in the form
\begin{equation}
\frac{d^2}{dt^2}\left(X(t)+c_2\rho(t)\right)\leq c_1X(t)+c_3\gamma\frac{d\rho (t)}{dt}+\frac{d^2F (t)}{dt^2}.
\end{equation}
From this inequality we obtain
\begin{eqnarray}
\label{X_ineq}
X(t)+c_2\rho(t)\leq  F(t)
 +c_1\int_{0}^{t}\int_{0}^{\sigma} X(\tau)\,d\tau d\sigma+c_3\gamma\int_0^t\rho(\tau)\,d\tau
\end{eqnarray}

In order to complete the proof,  we use the arguments of {\it reductio ad absurdum}. Let the conditions (\ref{cond1})--(\ref{cond2}) hold, but     a solution of the Cauchy problem  (\ref{Cauchy})   nevertheless exists for all  $t\in[0,T_0]$. Then  we can define
\begin{eqnarray}
T_1= \sup\left\{t\in[0,T_0] \colon X(s) \leq M(T_0) \mbox{\,\, for any \,\,} s\in [0, t]\right\}.
\end{eqnarray}
 From (\ref{X_ineq}) and (\ref{cond1})  for all $t\in [0, T_1]$  one has
\begin{eqnarray}
\label{X_ineq_1}
X(t)+c_2\rho(t)\leq F(t) + c_1  M(T_0) \frac{T_0^2}{2}
 +c_3\gamma\int_0^t\rho(\tau)d\tau
\nonumber \\
\label{X_ineq_2}
\leq M(T_0) -1  + c_1  M(T_0) \frac{T_0^2}{2}
 +c_3\gamma\int_0^t\rho(\tau)d\tau
\nonumber \\
< M(T_0)  + c_3\gamma\int_0^t\rho(\tau)d\tau.
\end{eqnarray}
Since $X(t) \geq 0$, by Gronwall's inequality
\begin{eqnarray}
\rho(t)\leq \frac{M(T_0)}{c_2}\exp\left(\frac{c_3}{c_2}\gamma
t\right).
\end{eqnarray}
Using this estimate back in right hand side of (\ref{X_ineq_1})  and using function $G(t)$ defined by (\ref{G}), we obtain
\begin{eqnarray}
\label{X_ineq_3}
X(t) \leq  F(t)+ G(T_0).
\end{eqnarray}
Therefore, using (\ref{cond2}) we conclude that $X(T_1) < M(T_0)$ and thus  $T_1=T_0$. Hence
\begin{eqnarray}
X(T_0)\leq F(T_0)+1<0,
\end{eqnarray}
which is impossible because $X(t)\geq 0$.

\section{``Early-collapse''}

Now we consider the case, when the NLS equation with gain, i.e. Eq. (\ref{NLS1}), contains  attractive nonlinearity, while the nonlinearity of the system with loss i.e. Eq.~(\ref{NLS2}), as well as nonlinear coupling, are either repulsive or zero,  i.e. we set
\begin{eqnarray}
\label{g_negative}
g_1>0, \quad g_2 \leq 0, \quad g\leq 0.
\end{eqnarray}

For this set of parameters, one can use  the idea of the ``early-collapse'' suggested in~\cite{KP}. The method is based on the fact that even if  the energy is a  growing function,  its growth  can   be controlled and it is hence  possible to choose the initial conditions
$u_0$ and $v_0$
such that the blowup occurs at sufficiently early times of the evolution.

In order to obtain the growth rate of the energy, we use (\ref{DE}) which subject to (\ref{g_negative}) allows one to obtain
\begin{eqnarray}
\label{E_estim_1}
  \frac{dE}{dt}\leq 2\gamma \int \left(
    |\nabla u|^2+|\nabla v|^2 +\kappa(\bar{u}v+u\bar{v}) -\frac{g_1}{2}|u|^4-\frac{g_2}{2}|v|^4-g|u|^2|v|^2\right) dx  +
    \nonumber
    \\
    +2\gamma\kappa (\|u\|_2^{{2}}+\|v\|_2^2) \leq 2\gamma E(t)+2\kappa\gamma S_0(t).
\end{eqnarray}
Thus using (\ref{S_estim}) and the definition of the energy (\ref{energy}) we derive
\begin{eqnarray}
\label{E_estim_2}
E(t)\leq
\left(E(0)+2\kappa\gamma S_0(0) t\right)e^{2\gamma t}=:E_{max}(t).
\end{eqnarray}
Next, we derive  the following estimate
\begin{eqnarray}
16\gamma\mbox{Im}\int\left(\bar{u}x\cdot\nabla u-\bar{v}x\cdot\nabla v\right)dx\leq
\frac{16\gamma^2}{N-2}X+4(N-2)  \left( \|\nabla u\|^2+\|\nabla v\|^2 \right),
\end{eqnarray}
and from the identity (\ref{X}) we  obtain
\begin{eqnarray}
\label{X_estim_new}
\frac{d^2 X}{dt^2}\leq c_4^2X+ 4N E_{max}(t)+4\kappa NS_0(t), \quad c_4=2\gamma\sqrt{\frac{\kappa}{\gamma}+\frac{N+2}{N-2}}.
\end{eqnarray}
Now one can obtain the upper bound for $X$:
\begin{eqnarray}
X(t)\leq Z(t)e^{c_4t},
\end{eqnarray}
where
\begin{eqnarray}
\label{X_estim_3}
Z(t):=X(0)+\int_{0}^{t}e^{-2c_4s}\left[Y(0)-c_4X(0)+4N\int_{0}^{s}e^{c_4\sigma}\left[E_{max}(\sigma)+\kappa S_0(0)e^{2\gamma \sigma}\right]d\sigma\right]ds,
\end{eqnarray}
and the function $E_{max}(\sigma)$  was defined in (\ref{E_estim_2}).

The obtained result can be  reformulated  as the following Theorem.
\begin{teo}
\label{theor2}
Let $N\geq 3$, the coefficients $g_1$, $g_2$ and $g$ satisfy (\ref{g_negative}),  and initial conditions  $u_0$ and $v_0$  in the Cauchy problem  (\ref{Cauchy})
are such that the function $Z(t)$ defined by (\ref{X_estim_3}) has a real positive zero $T_*$. Then  the solution of the problem (\ref{Cauchy}) does not exist in the interval $t\in [0, T_*]$.
\end{teo}

Since the occurrence of the blowup is now reduced to the study of the zeros of $Z(t)$ and this function depends on several parameters of the problem, we limit further consideration of this section by the numerical analysis. To this end now, as well as in all other numerical examples below, we limit the analysis to the Gaussian initial conditions (as the most typical for experimental settings):
\begin{eqnarray}
\label{IC_num}
u_0(x)=\frac{A}{\pi^{N/4} a^{N/2}}\exp\left(-\frac{|x|^2}{2a^2}\right),\quad
v_0(x)=\frac{B}{\pi^{N/4} b^{N/2}}\exp\left(-\frac{|x|^2}{2b^2}\right),
\end{eqnarray}
where $a$, $b$, $A$, and $B$ are positive constants.

\begin{figure} 
\centering
\includegraphics[width=0.3\textwidth]{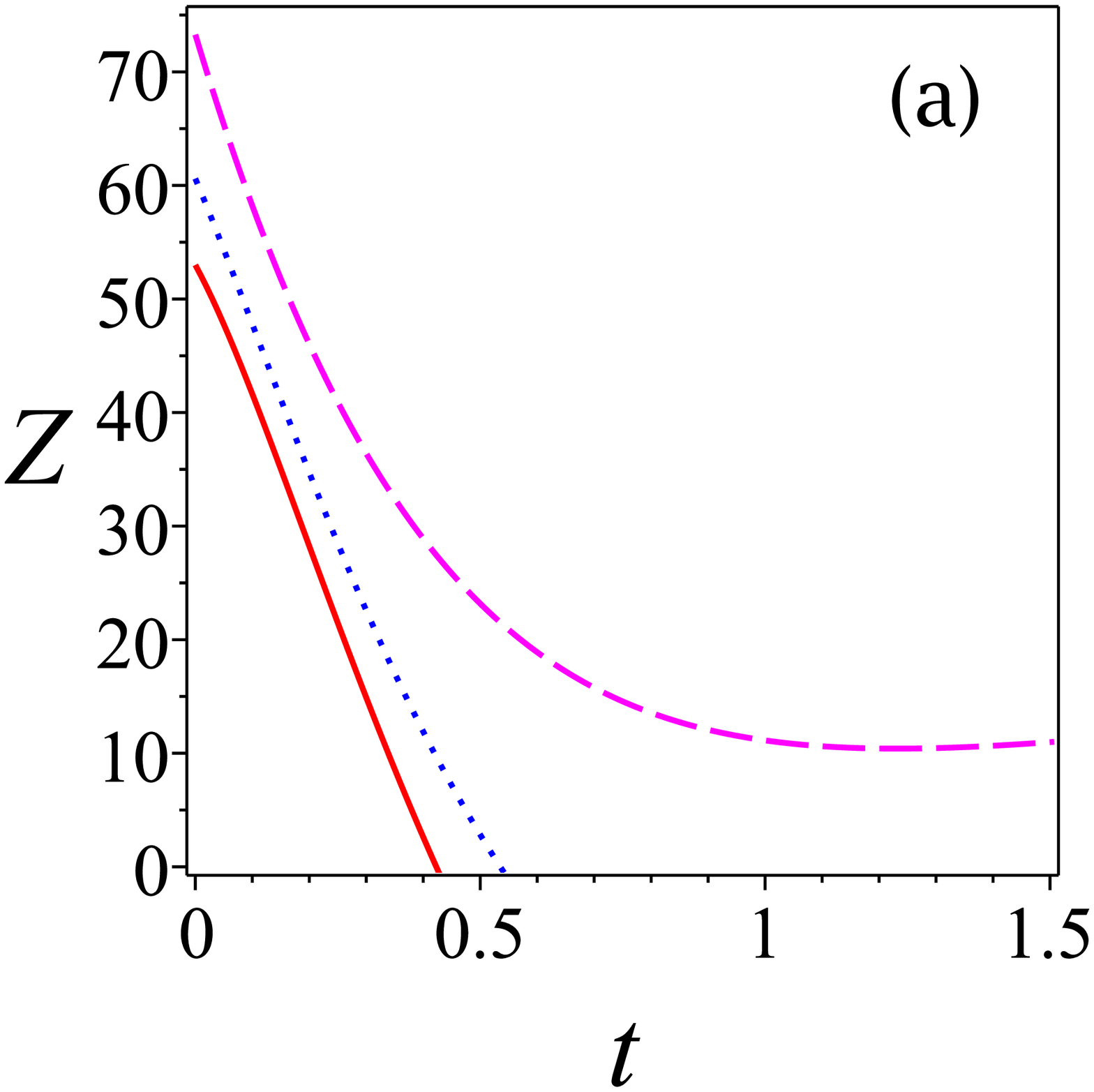}
\includegraphics[width=0.3\textwidth]{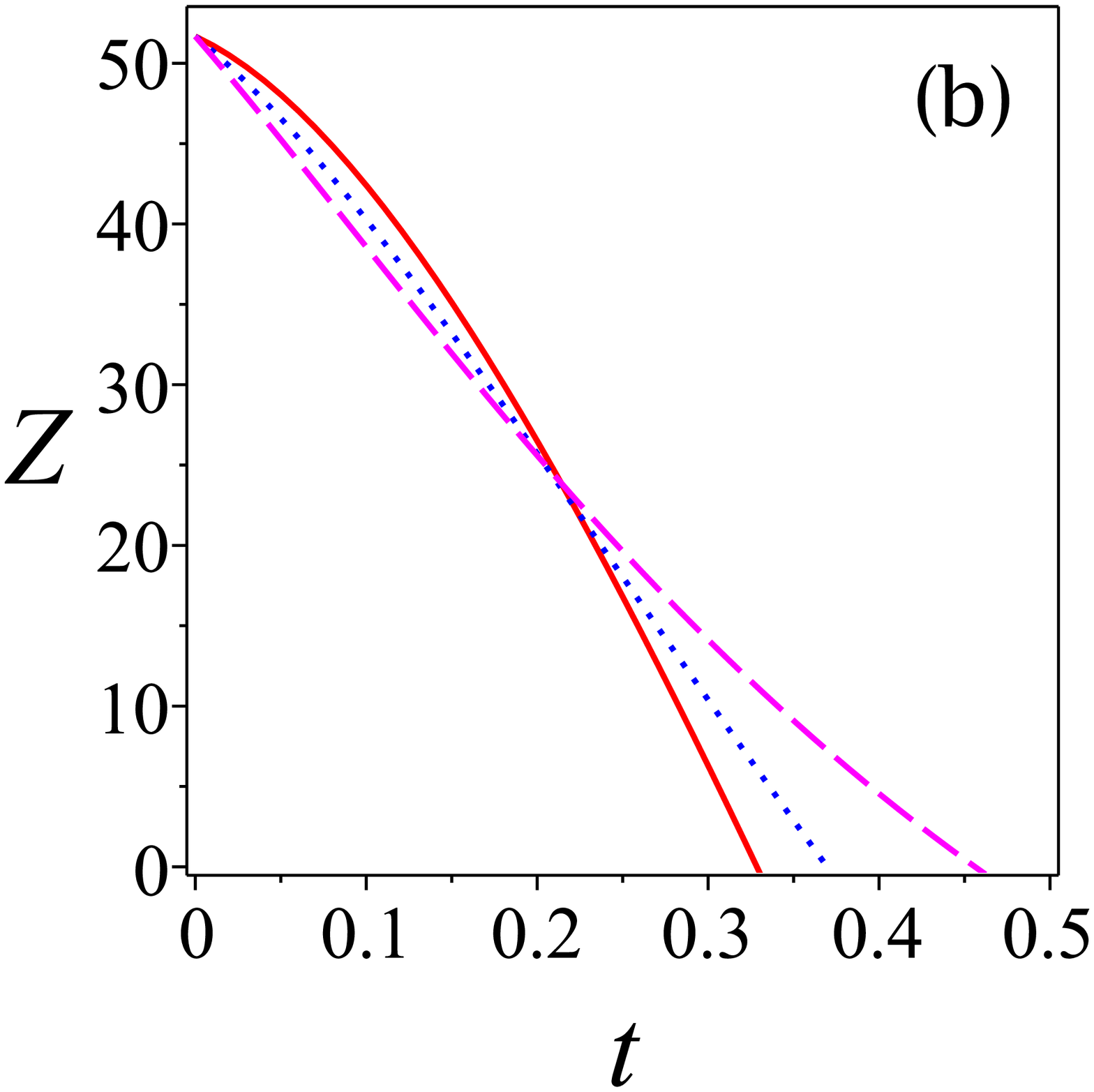}
\includegraphics[width=0.3\textwidth]{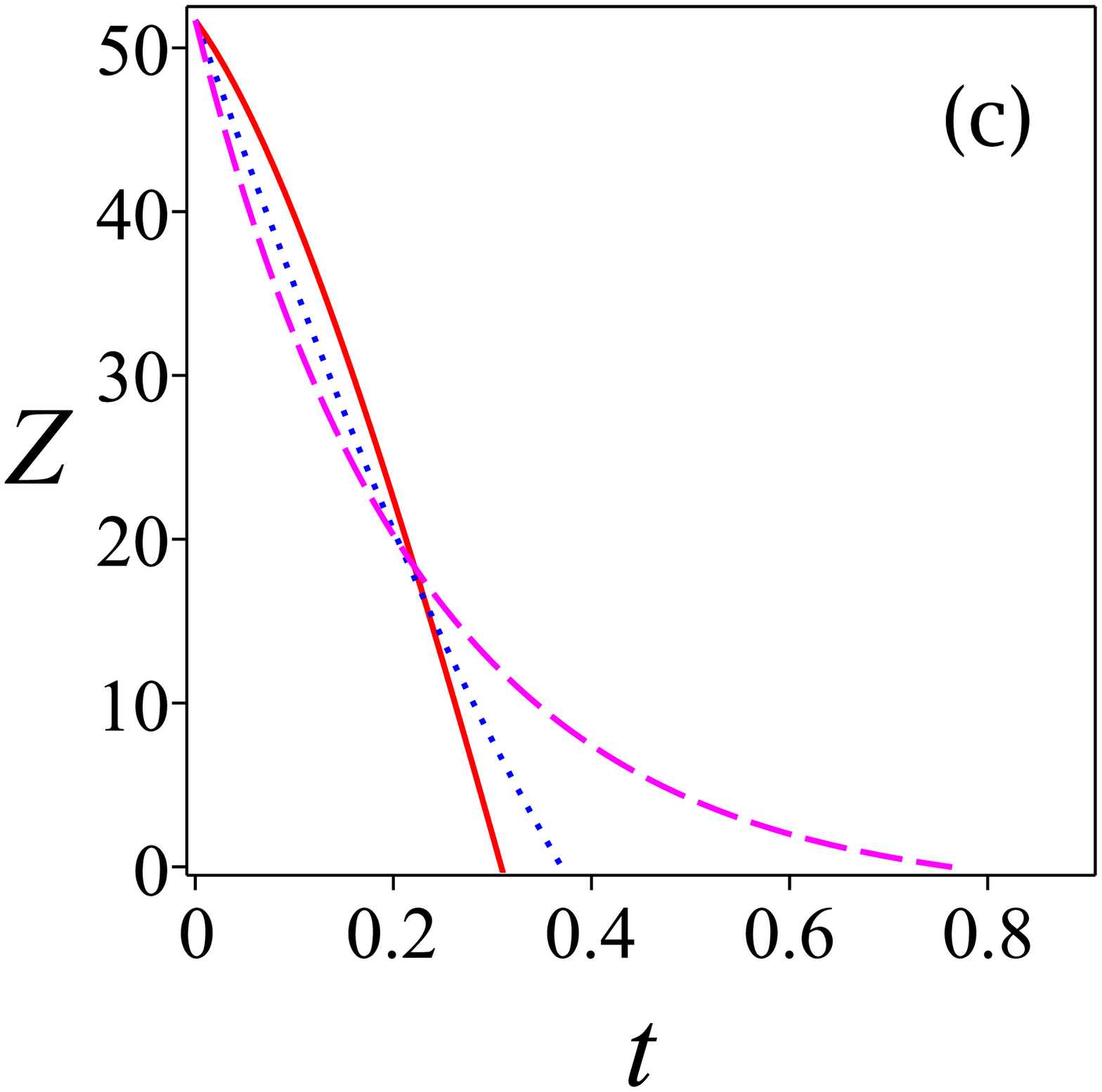}
\caption{Plots of $Z(t)$ for $N=3$, $g_1=4$, $g_2=-1$, $g=-0.5$, and the Gaussian initial conditions (\ref{IC_num}) with $A=5.8$ and $a=b=1$. The other parameters are as follows. Panel (a): $\gamma=0.5$, $\kappa=1$,
$B=1.3$ (solid line), $B=2.6$ (dotted line), $B=3.9$ (dashed line).
Panel (b)  $\gamma=0.5$,  $B=0.9$, $\gamma=0.15$ (solid line), $\gamma=0.3$ (dotted line), $\gamma=0.45$ (dashed line).
Panel (c): $\gamma=0.5$, $B=0.9$,  $\kappa=0.4$,   (solid line), $\kappa=0.8$ (dotted line), $\kappa=1.2$ (dashed line).
}
\label{fig:critical}
\end{figure}

In Fig.~\ref{fig:critical} we  show  plots of the  function  $Z(t)$ for several particular choices of the parameters in the case $N=3$, where we observe that for a proper choice of the initial conditions and the system parameters zeros of the function $Z(t)$ indeed exist. The time where $Z(t)$ becomes zero gives the respective $T_*$ (with is the minimal root obtained numerically).
In Fig.~\ref{fig:critical}(a) we show the behavior of $Z(t)$ for different values $B$, which corresponds to different initial conditions  $v_0(x)$ for the equation (\ref{NLS2}) with dissipation, provided the input of the NLS with gain, i.e. $u_0(x)$ is fixed. Increase of the initial amplitude $B$ of the pulse subjected dissipation [i.e   $v_0(x)$] results in larger values $T_*$ and eventually leads to nonexistence of the zeros of $Z(t)$.  In the cases where $T_*$ does not exist, the question about the finite time  blowup   remains open.

In Fig.~\ref{fig:critical}(b) we address different values of the  gain/loss coefficient $\gamma$, and in Fig.~\ref{fig:critical}(c) we illustrate different values of the coupling coefficient $\kappa$ at fixed initial conditions.
We observe increase of $T_*$ for stronger gain and dissipation [panel (b)], as well as for stronger coupling [panel (c)].

\section{A remark on the Manakov model}

The sufficient blowup conditions formulated in Sec.~\ref{sec:main}   allow for a simplification in the case of equal nonlinear coefficients:  $g_1=g_2=g>0$ (below we refer to this model as Manakov's one, after the original work~\cite{Manakov}).
Without loss of generality now we can set $g_1=g_2=g=1$ (which  can be  achieved by the renormalization $\sqrt{g} u\to u$, $\sqrt{g} v\to v$). 


An interesting property of the Manakov model is the existence of (at least) two integrals of motion. Indeed, straightforward computations show that the quantities
\begin{equation}
\label{int}
S_1(t)=S_1(0) \qquad \mbox{and} \qquad S(t):=\kappa S_0(t)-\gamma S_2(t)=\kappa S_0(0)-\gamma S_2(0)
\end{equation}
do not depend on time.
Then, from (\ref{difS0}) one obtains
\begin{eqnarray}
\frac{d^2S_0}{dt^2}+4 \omega^2 S_0= 4\kappa S,
\end{eqnarray}
where $\omega=\sqrt{\kappa^2-\gamma^2}$ as defined in (\ref{NLS1_rot}).
Thus in the case of unbroken the $\PT$-symmetric phase (\ref{unbroken}) the quantity $S_0(t)$ undergoes oscillatory dynamics:
\begin{eqnarray}
S_0(t)=\frac{\kappa S}{\omega^2}+S_{01}\cos(2\omega t)+S_{02}\sin(2\omega t), \\
S_{01}=S_0(0)\left(1-\frac{\kappa}{\omega^2}\right)+S_2(0)\frac{\gamma\kappa}{\omega^2}, \quad S_{02}=S_3(0)\frac{\gamma}{\omega},
\end{eqnarray}
and hence there exists   an upper bound for $S_0(t)$  [c.f. (\ref{S_estim})]:
\begin{equation}
 S_0(t)\leq S_{0,max} :=\frac{\kappa S}{\omega^2} + \sqrt{S_{01}^2+S_{02}^2}.
\end{equation}

In order to obtain sufficient conditions of the blowup in the Manakov case, we again employ   Eq.~(\ref{X}) and Eq.~(\ref{rearrange}) and,  using that the quantity $S_1(t)$ is conserved,  we can  replace the estimate (\ref{intermediate}) by the following one:
\begin{eqnarray}
\label{XMan}
\frac{d^2 X}{dt^2} + \frac{2(N-2)}{N+2}\left(2\|\nabla
u\|_2^2+2\|\nabla v\|_2^2 + Ng_1\|u\|_4^4+Ng_2\|v\|_4^4+ 2Ng \int
|u|^2|v|^2dx \right) \leq
\nonumber\\
\leq c_1X+\frac{16N}{N+2}(E(t) - \kappa S_1(0)).
\end{eqnarray}
One can repeat all the subsequent steps of the proof of Theorem~\ref{theor1} using the newly obtained estimate (\ref{XMan}). Therefore,
in the Manakov case  Theorem~\ref{theor1} can be reformulated with functions $F$, $M$ and $G$ replaced by $\hat{F}$, $\hat{M}$ and $\hat{G}$, respectively, where the new functions are defined as
\begin{eqnarray}
\label{FE}
\hat{F}(t) &=& X(0)+Y(0)t+\frac{8N}{N+2}(E(0) - \kappa S_1(0))t^2,
\\[3mm]
\hat{M}(t) &=& \sup_{\tau\in[0,t]}\hat{F}(\tau)+1,
\\[3mm]
\hat{G}(t)  &=& \hat{M}(t)\left(c_1 \frac{t^2}{2}+  \exp\left(\frac{48N\gamma t}{N+2}\right) - 1\right).
\end{eqnarray}

\section{Numerical illustrations}
\label{sec:numer}

The  analytical results obtained above give sufficient conditions for the finite time blowup, but do not describe
the  blowup dynamics and its dependence on the parameters  of the model, i.e.  $\gamma$, $\kappa$, $g_1$, $g_2$ and $g$.
In order to understand better the effect of those parameters on the phenomenon of blowup, now we resort to numerical analysis of  the Cauchy problem (\ref{Cauchy}).
We  concentrate on the case $N=3$ which is the most interesting one from the physical point of view. We also set $g_1=g_2$ and consider how the blowup depends on the nonlinear coupling $g$ and on the gain/loss coefficient $\gamma$. As we mentioned above, without loss of generality we can set $g_1=g_2=1$ (which is achieved   by the renormalization $\sqrt{g_1} u\to u$, $\sqrt{g_2} v\to v$,  $g/g_1=g/g_2 \to g$). Using another evident  renormalization, without loss of generality we can set  $\kappa=1$. Therefore we concentrate our attention on the system
\begin{subequations}
\label{NLS_3D}
\begin{eqnarray}
\label{NLS1_3D}
&i u_t=-\Delta u +i\gamma u+ v-(|u|^2+g|v|^2)u, &
\\
\label{NLS2_3D}
&i v_t=-\Delta v -i\gamma v+ u-(g|u|^2+|v|^2)v,  &
\end{eqnarray}
\end{subequations}
 where $x\in\mathbb{R}^3$ and $ t\geq 0$.

\begin{figure}
\centering
{\includegraphics[width=\textwidth]{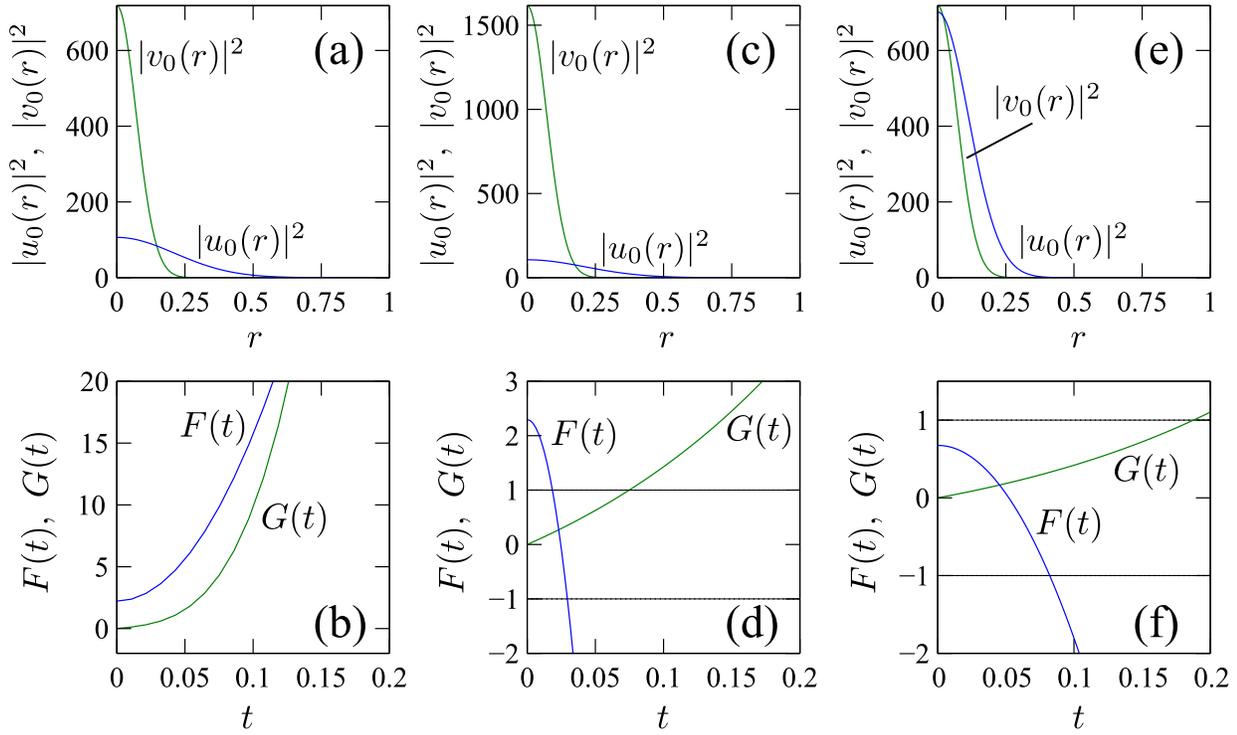}}
\caption{(a) and (b) Plots of initial conditions (\ref{IC_num}) for $N=3$,  $A=4$, $B=2$, $a=0.3$  and $b=0.1$ (with $r=|x|$) and behavior of  the associated functions $F(t)$ and $G(t)$ defined by (\ref{F}) and (\ref{G}) with $\kappa=1$, $\gamma=0.5$,  $g_1=g_2=1$ and $g = -0.5$; (c) and (d) show the same as (a) and (b) but with $B=3$; (e) and (f) show the same as (a) and (b) but with $b=0.16$. Notice that the blowup conditions (\ref{cond1})--(\ref{cond2}) are not satisfied in (a)--(b), but become satisfied in (c)--(d) and (e)--(f).
}
\label{fig-IC}
\end{figure}

We first notice that  the blowup conditions (\ref{cond1})--(\ref{cond2}) can be easily satisfied by a proper choice of  the  initial pulses $u_0$ and $v_0$. Several examples are presented in Fig.~\ref{fig-IC} which shows spatial profiles of initial conditions $u_0(x)$ and $v_0(x)$ chosen in the form of the Gaussian beams (\ref{IC_num}) and the behavior of the associated functions $F(t)$ and $G(t)$ defined by (\ref{F}) and (\ref{G}). Panels  Fig.~\ref{fig-IC}~(a--b)  correspond to a situation when the blowup conditions are not satisfied as functions $F(t)$ and $G(t)$ grow monotonously. However, increasing   powers $A$ and $B$ of the input beams [or one of the beams, see  Fig.~\ref{fig-IC}~(c--d)] or decreasing the characteristic widths $a$ and $b$ of the beams [Fig.~\ref{fig-IC}(e--f)], one can easily satisfy the blowup conditions.

We also performed numerical simulations of the three-dimensional system  (\ref{NLS_3D}) in the spherical symmetric case, i.e. assuming that the functions $u$ and $v$  depend only on the radius $r=|x|$ in the spherical coordinates and do not depend on the polar and azimuthal angles. In this case it is convenient to introduce    new functions $p(r, t)=ru$ and $q(r, t)=rv$ and to reformulate the problem as follows
\begin{subequations}
\label{shp}
\begin{eqnarray}
\label{shp1}
ip_t &=& -p_{rr} + i\gamma p + q - r^{-2}(|p|^2+g|q|^2)p,\\
\label{shp2}
iq_t &=& -q_{rr} - i\gamma q + p - r^{-2}(g|p|^2+|q|^2)q.
\end{eqnarray}
\end{subequations}
We also   limited the  numerical  study to a finite interval $r\in[0,L]$ subject to the zero boundary conditions $p(0, t)=q(0,t)=p(L, t)=q(L,t)=0$. The interval width   $L$ is taken sufficiently large  such that  increase of $L$ practically does not affect the results presented below.

We solved   system  (\ref{shp}) using a semi-implicit finite-difference scheme proposed in \cite{Trofimov}. The numerical simulations  were carried on up to the time $t_*$ when the ratios $|u(t_*,0)|/|u_0(0)|$ and $|v(t_*,0)|/|v_0(0)|$ were of order $10^2$, and then the simulations were interrupted. We used the spatial step of order $10^{-5}$ and the adaptively decreasing temporal step which was of order $10^{-5}$ in the beginning of simulations (i.e. at $t=0$) and of  order   $10^{-7}$ at the time $t=t_*$ of the termination of the simulations. We also performed several additional runs with smaller spatial steps and checked convergence of the numerical solution.

In the numerical simulations we have observed different scenarios of the dynamics: dispersion of the initial pulses, growth  of the solution  and its  derivative  at the origin $x=0$   in one of the components (either with gain or with dissipation), and the simultaneous growth in  both the components. The dispersion corresponds to the spreading of the initial pulses (which eventually occupy the entire computational domain $[0, L]$), while the observed growth is a presumable  manifestation of the blowup.
Our numerical results allow us to conjecture that depending on the initial conditions and parameters of the model the blowup can occur either in both or in only one component. In the latter case this can be a component either with gain or with damping.

Summarizing results of dynamical simulation of   system, we first of all   conjecture    that the obtained above blowup conditions  (\ref{cond1})--(\ref{cond2})   are not necessary for the  blowup to occur: it is easy to find  initial conditions that do not satisfy (\ref{cond1})--(\ref{cond2}) but nevertheless   grow  in numerical simulations.  Moreover, our numerics indicate that the growth  can occur even for initial conditions with $E(0) \geq 0$ and $Y(0)\geq 0$, i.e. in the case when conditions  (\ref{cond1})--(\ref{cond2}) obviously can not be satisfied for any $T_0$.

\begin{figure}
\centering
{\includegraphics[width=\textwidth]{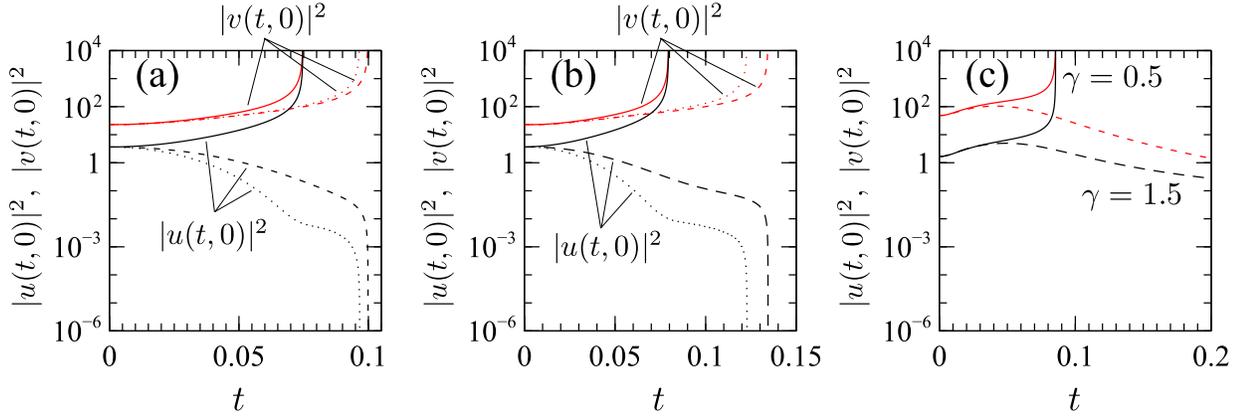}}
\caption{(a)-(b) Numerically obtained dynamics  with the initial
conditions (\ref{IC_num}) with $A=4.5$, $B=4$,  $a=1$ and $b=0.5$
for $\gamma=0.5$  (a) and $\gamma=1.5$ (b). Solid, dashed and
dotted lines correspond to $g=1$, $g=-1$ and $g=-2$, respectively.
(c) The dynamics obtained for $A=0.5$, $B=2.7$, $a=b=0.3$ and
$g=1$. The solution blows up for the unbroken $\PT$
symmetry ($\gamma=0.5$, solid curves) but disperses (at least, at the initial stage of the evolution) for the broken
$\PT$ symmetry ($\gamma=1.5$, dashed lines). The dependencies are
shown in the semi-logarithmic scale. } \label{fig-1}
\end{figure}

Several (the most interesting) examples of the blowup-like dynamics obtained from numerical solution of   system (\ref{shp}) are presented in Fig.~\ref{fig-1}. The case of the unbroken $\PT$ symmetry is illustrated in Fig.~\ref{fig-1}(a) from where one observes that positive and negative values of $g$ result in different scenarios of the blowup-like dynamics. Namely, for $g = 1$ (this is the attractive inter-species nonlinearity) we observe that the   the intensity growth  occurs in both the components. However for  $g=-1$ and $g=-2$ the same initial conditions display  growth only in the \textit{dissipative} component $v(x, t)$ while in  the gain component $u(x,t)$ the intensity at the origin decays. Notice also that the maximal squared amplitude in the component with gain  $\max_{r\in [0, L]} |u(x,t)|^2$   for $g=-1$ and $g=-2$ (not shown on the plot) was of the order of $1$  at the moment when the simulations were terminated. It is also interesting to notice that for $g=-1$ and $g=-2$ the only component ``responsible'' for the blowup is the one with dissipation  [i.e. $v(x,t)$]  and that the growth   with $g=-2$ (the dotted curves) occurs earlier than that with $g=-1$ (the dashed curves) in spite of the additional ``defocusing''  that large negative values of  $g$ induce. We finally notice that the case $g=-2$ is not covered by the condition of the Theorem~1, which is another illustration for the fact   that the obtained blowup conditions (\ref{cond1})--(\ref{cond2}) are sufficient but not necessary.

Repeating the same simulations for $\gamma=1.5$ [the broken $\PT$ symmetry, Fig.~\ref{fig-1}(b)] we do not observe  any quantitative difference in the blowup scenarios with respect to the case $\gamma=0.5$. However it is interesting to notice that for the chosen initial conditions  and for the three  different values of $g$ addressed in Fig.~\ref{fig-1}(b) the numerical blowup for the broken $\PT$ symmetry occurs \textit{later} than that with the unbroken $\PT$ symmetry in Fig.~\ref{fig-1}(a) (this is especially well-visible for $g=-1$ and $g=-2$). However this behavior is conditioned by the chosen initial profiles. It is also possible to find  initial conditions   for which the broken $\PT$ symmetry precipitates  the intensity growth  (with respect to the unbroken $\PT$ symmetry). Another interesting effect is that it is  possible to select the initial conditions which grow   for the unbroken $\PT$ symmetry but disperse (at least, at the initial stage of the evolution) for the broken $\PT$ symmetry [Fig.~\ref{fig-1}(c)],  which is an   indication of absence of any direct relation between the blowup phenomenon and the $\PT$ symmetry breaking.

\begin{figure}
\centering
{\includegraphics[width=0.67\textwidth]{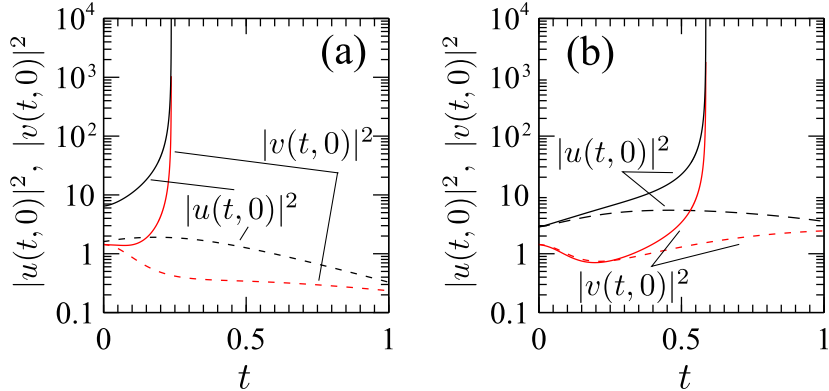}}
\caption{Change of the dynamics after modification of parameter $A$ of the initial conditions [panel (a)] and parameter of $\PT$ symmetry $\gamma$ [panel (b)].
 The dynamics in (a) is obtained for  the initial conditions (\ref{IC_num}) with  $B=1$,  $a=1$ and $b=0.5$,
   $A=3$ (dashed curves) and $A=6$ (solid curves), and with  $\gamma=0.5$.
 The dynamics in (b) is obtained for the initial
conditions (\ref{IC_num}) with  $A= B=1$,  $a=1$, $b=0.5$, and with $\gamma=0.5$
     (dashed curves) and $\gamma=0.9$ (solid curves).  In both panels $g=1$. The dependencies are
shown in the semi-logarithmic scale. } \label{fig-referee}
\end{figure}

In Fig.~\ref{fig-referee} we demonstrate how increase of the amplitude $A$ of the initial conditions  [Fig.~\ref{fig-referee}(a)] and increase of  the parameter of the $\PT$ symmetry $\gamma$  [Fig.~\ref{fig-referee}(b)] changes the dynamics from the dispersive (and presumably global) one to the blowup.

\section{Concluding remarks}

In this paper we have established sufficient conditions for finite-time blow up of initial pulses whose evolution is described by the coupled $\PT$-symmetric NLS equations, one of which has linear gain and another one has dissipation. The equations include both linear and nonlinear coupling. The proofs were presented for the two cases.  In the first situation the both NLS equations have attractive (focusing) self-phase nonlinearities and either attractive or weakly repulsive cross-phase nonlinearities (i.e. nonlinear coupling coefficients). In the second case the equation with gain was self-attracting while all other nonlinear interactions were repulsive (defocusing).

Neither proof nor the results rely on whether the $\PT$ symmetric phase is broken or not. In the formulations of theorems the gain/dissipation strength and the coupling enter as independent parameters (unlike this happens in the linear problem, where their relation is the only significant parameter). This is however expectable in view of the fact that the blowup is an essentially nonlinear phenomenon.

We also performed numerical studies of the initial stages of evolution of the initial data.
The consideration was limited to the simplest Gaussian shapes; comprehensive numerical study of the initial data is left as an open question. However even these first studies revealed interesting counter-intuitive dynamics like dispersion of pulses in the case of broken (linear) symmetry and blow up of the same pules in the case of unbroken symmetry, or occurrence of collapse of the field in the equation with dissipation with the field subjected to gain reaming bounded.

Finally, we notice that the developed theory does not allow for direct generalization for the critical collapse, which is left, so far, as another open problem.


\section*{ Acknowledgments}
JPD and MF acknowledge support of the FCT (Portugal) grant PEst-OE/MAT/UI0209/2013. VVK and DAZ acknowledge  support of the FCT (Portugal) grants PEst-OE/FIS/UI0618/2014 and   PTDC/FIS-OPT/1918/2012.

\end{document}